\documentclass{article}
\usepackage[ngerman,french,english]{babel}
\usepackage{amsmath,amssymb,amsthm,amscd, mathtools, latexsym, tikz-cd}
\usepackage{enumerate,varioref, dsfont}           
\usepackage{enumitem}

\usepackage{hyperref}
\usepackage{url}

\usepackage{hyperref}
\usepackage{url}

\newtheorem{Thm}{Theorem}[section]
\newtheorem{Lem}{Lemma}[section]

\newtheorem{Cor}{Corollary}[section]
\newtheorem{Conj}{Conjecture}[section]

\newtheorem{claim}{Claim}[section]

{\theoremstyle{plain}

\newtheorem{Not}{Notation}[section]
}
\renewcommand{\ker}{\text{ker}}

\def\cit{{\mathbb C}}

\def\qit{{\mathbb Q}}

\def\pit{{\mathbb P}}

\def\0{{\mathcal O}}

\def\E{{\mathcal E}}

\begin{document}
\title{The Standard Conjectures for the Variety of Lines on a Cubic Hypersurface}

\author{HUMBERTO DIAZ}

\newcommand{\Addresses}{{\bigskip \footnotesize
\textsc{Department of Mathematics, University of California-Riverside, Riverside, CA 92521} \par \nopagebreak
\textit{Email address}: \ \texttt{humbertd@ucr.edu}}}

\date{}

\maketitle

\begin{abstract}
The purpose of this note is to prove Grothendieck's standard conjectures for the Fano variety of lines on a smooth cubic hypersurface in projective space. 
\end{abstract}

\maketitle

\section*{Introduction}
\noindent Let $X \subset \pit^{n+1}$ be a smooth cubic hypersurface over a field $k$ (with $n \geq 3$). The Fano variety of lines on $X$, $F:=F(X)$, has proved to be useful in understanding the geometry of the cubic. A thorough study of the Fano variety was completed by Altman and Kleiman in \cite{AK}, in which they show, among other things, that the Fano variety of $X$ is smooth and has dimension $2(n-2)$. In the case of $n=3$, the Fano variety is a surface of general type which possesses a great many remarkable properties, which were used by Clemens and Griffiths (in their well-known paper \cite{CG}) to establish that the smooth cubic threefold is not rational (when $k = \cit$). For $n=4$, Beauville and Donagi showed in \cite{BD} that $F$ has the structure of a hyper-K\"ahler variety by establishing that it is deformation equivalent to $S^{[2]}$, the second punctual Hilbert scheme of a $K3$ surface. Quite recently, Galkin and Schinder in \cite{GS} were able to prove the following relation in the Grothendieck ring of varieties:
\begin{equation} [X^{[2]}] = [X][\pit^{n+1}] + \mathds{L}^{2}[F] \in K_{0} (Var)\label{decomp} \end{equation}
where $X^{[2]}$ denotes the second punctual Hilbert scheme of $X$. By applying the Hodge realization map $K_{0} (Var) \to K_{0} (\qit$-HS) to (\ref{decomp}), they also obtain an (abstract) isomorphism of $\qit$-Hodge structures:
\begin{equation} H^{*} (F, \qit) \cong \text{Sym}^{2} H^{n}_{prim} (X, \qit) (2) \oplus \bigoplus_{k=0}^{n-2} H^{n}_{prim} (X, \qit) (-k+1) \oplus \bigoplus_{k=0}^{2(n-2)} \qit (-k)^{\oplus a_{\frac{k}{2}}}\label{cohom-1} \end{equation}
where $H^{n}_{prim} (X, \qit)$ is the primitive cohomology and $a_{\frac{k}{2}}$ are some positive integers that will be given in the sequel.\\
\indent Using decomposition (\ref{cohom-1}) as a guide will allow us to prove the standard conjectures for $F$. These conjectures were first stated by Grothendieck in his paper \cite{Gr} and concern the existence of certain fundamental algebraic cycles. 

\begin{Conj}\label{Gr}(Grothendieck) Let $Y$ be a smooth projective variety of dimension $d_{Y}$ over $k$, $H^{*}$ denote a Weil cohomology and $H^{2j}_{alg}(Y) \subset H^{2j} (Y)$ denote the image of the cycle class map, $CH^{j} (Y) \to H^{2j}(Y)$.
\begin{enumerate}
\item\label{B}(Lefschetz standard) For an ample divisor $h \in Pic(Y)$, there exist correspondences \[\Lambda_{j} \in CH^{j}(Y \times Y)\] such that $\Lambda_{j*}: H^{2d_{Y}-j}(Y) \to H^{j}(Y)$ is the inverse of the Lefschetz isomorphism:
\[ L_{h}^{j} := \cdot h^{j}: H^{j} (Y) \to H^{2d_{Y}-j} (Y)\]
is induced by an algebraic cycle.
\item\label{C}(K\"unneth) The K\"unneth components $\delta_{j} \in H^{j} (Y) \otimes H^{2d_{Y}-j} (Y)$ in the decomposition of the diagonal:
\[ \Delta_{Y} \in H^{2d_{Y}} (Y \times Y) \cong \bigoplus_{j=0}^{2d_{Y}} H^{j} (Y) \otimes H^{2d_{Y}-j} (Y) \]
are induced by algebraic cycles; i.e., $\delta_{j} \in H^{2d_{Y}}_{alg} (Y \times Y)$.
\item\label{D}(${\sim}_{num} = {\sim}_{hom}$) For any $\gamma \in H^{2j}_{alg} (Y)$, if $\gamma\cdot \gamma' = 0$ for all $\gamma' \in H^{2(d_{Y}-j)}_{alg} (Y)$, then $\gamma =0$.
\end{enumerate}
\end{Conj}

\noindent The conjectures are known in a few important cases, but in general remain wide open. They are known to be true, for instance, when $X$ is an Abelian variety (\cite{Kl}), as well as when the cycle class map is an isomorphism (\cite{Kl2}). (This condition holds for varieties admitting a cellular decomposition such as Grassmannians, toric varieties, etc.) It is also known that the conjectures remain true for hyperplane sections, products, and projective bundles. It is noteworthy that when $k$ is of characteristic $0$, Conjecture \ref{Gr} \ref{B} implies the other 2 conjectures (see \cite{Kl} for details). Also, Arapura proves that the standard conjectures for uniruled threefolds, unirational fourfolds, the moduli space of stable vector bundles over a smooth projective curve, and for the Hilbert scheme $S^{[n]}$ of every smooth projective surface (\cite{A}, Corollaries 4.3, 7.2 and 7.5). Finally, Charles and Markman have shown in \cite{CM} that the conjectures are true for all smooth projective varieties which are deformation equivalent to a Hilbert scheme of $K3$ surfaces.
\begin{Thm}\label{main} The standard conjectures hold for the Fano variety of lines of a smooth cubic hypersurface in $\pit^{n+1}_{\cit}$ when $H^{*}$ is singular cohomology with $\qit$ coefficients.
\end{Thm}

\noindent The idea will be to show that the direct sum of Hodge classes $\qit(-k)^{a_{\frac{k}{2}}}$ in (\ref{cohom-1}) is algebraic and that the orthogonal complement of this sum is an appropriate Tate twist of $H^{n}_{prim} (X)$. This will allow us to construct correspondences $\Gamma_{k} \in CH^{k} (F \times F)$ which induce isomorphisms:
\[ \Gamma_{k*}: H^{4(n-2)-k} (F) \to H^{k} (F) \]
for $k \leq 2(n-2)$. Then, arguments of Kleiman from \cite{Kl} can be used to deduce the standard conjectures. 
\subsection*{Conventions} Throughout this note, $X$ will denote a smooth cubic hypersurface over $\cit$ and $F$ its Fano variety of lines, $G_{1,n+1}$ the Grassmannian of lines in $\pit^{n+1}$ and $H = c_{1} (\mathcal{O}_{X}(1)) \in Pic (X)$. We will let $H^{*}$ denote singular cohomology with  coefficients $\qit$ coefficients. (However, most of the results hold over any field with any choice of Weil cohomology.) We will also let $CH^{*}$ denote the Chow group with $\qit$ coefficients and for a correspondence $\phi \in CH^{*} (Z \times Y)$, we will let 
$\phi_{*}: H^{*} (Z) \to H^{*} (Y)$ denote the usual action on cohomology. For convenience, we will abuse notation by omitting Tate twists. We will also let $[x]$ denote the greatest integer $\leq x$. 

\section{Lemmas} 
\noindent In this section, we will prove a few geometric facts about the cubic hypersurface and its Fano variety of lines, which will be necessary for the proof of Theorem \ref{main}.\\

\noindent The universal line over the Fano variety $F$ gives a projective bundle:
\[ p: \pit(\E) \to F \]
and, hence, a natural imbedding $\iota: F \hookrightarrow G_{1,n+1}$ into the Grassmannian of lines on $\pit^{n+1}$. We can use the Grassmannian to account for the Hodge classes which appear in (\ref{cohom-1}); more precisely, we have the following:
\begin{Lem}\label{first-lem} Let $i: Y \hookrightarrow G_{1,n+1}$ be smooth closed subvariety of the Grassmannian of codimension $r$. Then, for $k \leq 2(n+2)-2r$, the pull-back \[\iota^{*}: H^{k} (G_{1,n+1}) \to H^{k} (Y)\] is injective. 
\begin{proof} Note that $H^{*} (G_{1,n+1})$ is (as a graded algebra) generated by $\{ c_{1} (\E), c_{2} (\E) \}$. Since the class of $Y$ is effective in $G_{1,n+1}$, it represents a (non-zero) homogeneous degree $2r$ polynomial in the Chern classes. So, suppose that $\alpha \in H^{k} (G_{1,n+1})$ is a (non-zero) polynomial in the Chern classes such that $\iota^{*}\alpha = 0 \in H^{k} (Y)$. Then, from the projection formula, we deduce that
\[ 0 = \iota_{*}\iota^{*}\alpha = \alpha \cdot Y \in H^{k+2r} (Gr(1, n+1)) \]
This means that $\alpha \cdot Y$ is a non-zero polynomial in the Chern classes of degree $=k+2r \leq 2(n+2)$. But this is already a contradiction thanks to the following claim:
\begin{claim} There are no relations among the Chern classes in degrees $\leq 2(n+2)$. 
\begin{proof}[Proof of Claim] For $2d \leq 2(n+2)$, the polynomials in the Chern classes are generated by the set 
\[\{ c_{1}(\E)^{d}, c_{2}(\E)\cdot c_{1}(\E)^{d-2}, \ldots , c_{2}(\E)^{[\frac{d}{2}]}\cdot c_{1}(\E)^{d-2[\frac{d}{2}]} \}\]
We see then that the rank of $H^{2d} (G_{1,n+1})$ is $\leq [\frac{d+2}{2}]$. However, by counting Schubert classes, one deduces that the rank is given by the number of solutions to:
\[ n_{1} + 2n_{2} = d \]
where $n_{1}, n_{2}$ are non-negative integers satisfying $n_{1} + n_{2} \leq n+2$. This latter is no restriction so long as $d \leq n+2$. It follows that the rank is exactly $[\frac{d+2}{2}]$. Thus, the claim.
\end{proof} 
\end{claim} 
\end{proof}
\end{Lem}
\begin{Cor}\label{first} For $k \leq 2(n-2)$, the pull-back \[\iota^{*}: H^{k} (G_{1,n+1}) \to H^{k} (F)\] is injective. Moreover, for $0 \leq k \leq 4(n-2)$ with $k \neq 2(n-2)$,
\[ \text{rank}(\iota^{*}H^{k} (G_{1,n+1})) = \left\{ \begin{array}{cc} \geq a_{\frac{k}{2}} & \text{ for $k$ even}\\ 0 & \text{ for $k$ odd} \end{array} \right.\]
\begin{proof} The first statement follows from Lemma \ref{first-lem} and the fact that $F$ has codimension $4$ in $G_{1,n+1}$. The second statement is trivially true when $k$ is odd. On the other hand, when $k$ is even, we note that by Corollary 5.7 of \cite{GS}, we have 
$$a_{\frac{k}{2}} = \left\{\begin{array}{cc}  \text{$[\frac{k+4}{4}]$}  & \text{ for $k < 2(n-2)$}\\ \text{$[\frac{4(n-2)-k+4}{4}]$} & \text{ for $k > 2(n-2)$}\end{array} \right.$$
and we observe that $a_{\frac{k}{2}} = a_{\frac{4(n-2)-k}{2}}$. Now, the proof of the claim above shows that for $k < 2(n-2)$ we have
\[ \text{rank}(\iota^{*}H^{k} (G_{1,n+1})) = [\frac{k+4}{4}] = a_{\frac{k}{2}} \] 
\noindent On the other hand, for $k > 2(n-2)$ then by the Hard Lefschetz theorem, we have 
\[\begin{split} \text{rank}(\iota^{*}H^{k} (G_{1,n+1})) &\geq \text{rank}(\iota^{*}H^{4(n-2)-k} (G_{1,n+1})\cdot c_{1}(\E)^{k-2(n-2)})\\ & = \text{rank}(\iota^{*}H^{4(n-2)-k} (G_{1,n+1}))\\ & = a_{\frac{4(n-2)-k}{2}} = a_{\frac{k}{2}} \end{split} \]
\end{proof}
\end{Cor}

\begin{Cor}\label{cor-cor} Whenever $k<n-2$ or $k>3(n-2)$ or $n-2 \leq k \leq 3(n-2)$, $k\neq 2(n-2)$ and $k \not\equiv n \pmod 2$,
\begin{equation} H^{k} (F) = \iota^{*}H^{k} (G_{1,n+1})\label{Tate} \end{equation}
\begin{proof} In all these cases, an inspection of (\ref{cohom-1}) shows that 
\[ \text{rank}(H^{k} (F)) = \left\{ \begin{array}{cc}  a_{\frac{k}{2}} & \text{ for $k$ even}\\ 0 & \text{ for $k$ odd} \end{array} \right.\]
which gives $\text{rank}(H^{k} (F)) \leq \text{rank}(\iota^{*}H^{k} (G_{1,n+1}))$, according to Corollary \ref{first}. So, this must in fact be an equality and this gives (\ref{Tate}).
\end{proof}
\end{Cor}

\noindent There are also natural projection maps:
\begin{equation}
\begin{tikzcd}
\pit(\E) \arrow{r}{p} \arrow{d}{q} & F\\
X
\end{tikzcd}
\end{equation}
This gives the well-known {\em cylinder correspondence} $\Gamma \in CH^{n-1} (X \times F)$ whose action on cohomology is given by $\Gamma_{*} = p_{*}q^{*}$, as well as its transpose $\prescript{t}{}{\Gamma}_{*} = q_{*}p^{*}$. Moreover, we have the following result:
\begin{Lem}\label{second} Let $H^{n}_{prim} (X)$ denote the primitive cohomology of $X$. Then, \[\Gamma_{*} = p_{*}q^{*}: H^{n}_{prim} (X) \to H^{n-2} (F)\] is injective. When $n$ is odd, this is an isomorphism. When $n$ is even, ${\Gamma}_{*}H^{n}_{prim} (X)$ and $\iota^{*}H^{3(n-2)} (G_{1,n+1})$ are orthogonal with respect to the cup product.
\begin{proof} The first statement is perhaps well-known. However, since the author cound not find a reference, we proceed as in \cite{BD}. Indeed, since $p: \pit(\E) \to F$ is a $\pit^{1}$-bundle, there is a decomposition:
\[ H^{n} (\pit(\E)) = p^{*}H^{n} (F) \oplus  p^{*}H^{n-2} (F)\cdot h\]
where $h = c_{1} (\mathcal{O}(1))$ is the class of the (anti-)tautological bundle. 
From \cite{F} Chapter 3 (or otherwise), we have
\[ p_{*}(p^{*}\gamma_{n}) = 0 , \ p_{*}(p^{*}\gamma_{n-2}\cdot h) = \gamma_{n-2} \]
for $\gamma_{j} \in H^{j} (F)$. Since $q^{*}$ is injective, it suffices to show that
\[ p^{*}H^{n} (F) \cap  q^{*}H^{n}_{prim} (X) = 0 \]
So, suppose there is some $\gamma \in H^{n}_{prim} (X)$ and some $\gamma_{n} \in H^{n} (F)$ such that
\[ p^{*}\gamma_{n} = q^{*}\gamma \in H^{n} (\pit(\E)) \]
We then observe that $h = q^{*}H$, from which it follows that
\begin{equation}
p^{*}\gamma_{n-2}\cdot h = q^{*}\gamma\cdot q^{*}H = q^{*} (\gamma \cdot H)
\label{key} \end{equation}
Since $\gamma \in H^{n}_{prim} (X)$, $\gamma \cdot H =0 \in H^{n+2} (X)$ and then (\ref{key}) gives
\[ \gamma_{n-2} = p_{*}(p^{*}\gamma_{n-2}\cdot h) = p_{*}(q^{*}(\gamma \cdot H)) = 0 \]
This gives the first statement. When $n$ is odd, an inspection of the decomposition in (\ref{cohom-1}) reveals that $H^{n-2}(F)$ has the same dimension as $H^{n}(X)$, which implies that $\Gamma_{*}$ is an isomorphism in this case. When $n$ is even, we would like to prove:
\begin{equation}\iota^{*}H^{3(n-2)} (G_{1,n+1}) \cdot {\Gamma}_{*}H^{n}_{prim} (X) = 0\label{goal0} \end{equation}
To this end, we have the following commutative diagrams:
\begin{equation}\begin{CD}
F  @<{p}<< \pit(\E) @>{q}>> X\\
@V{\iota}VV @V{\overline{\iota}}VV @V{j}VV\\
G_{1,n+1}  @<{\overline{p}}<< \pit_{G_{1,n+1}}(\E) @>{\overline{q}}>> \pit^{n+1}\label{comm}
\end{CD}\end{equation} 
where $\pit_{G_{1,n+1}}(\E) \xrightarrow{\overline{p}} G_{1,n+1}$ denotes the tautological $\pit^{1}$-bundle and $j: X \hookrightarrow \pit^{n+1}$ is the inclusion. Then, (\ref{goal0}) becomes 
\[\begin{split} \iota^{*}H^{3(n-2)} (G_{1,n+1}) \cdot {\Gamma}_{*}H^{n}_{prim} (X) & = \iota^{*}H^{3(n-2)} (G_{1,n+1}) \cdot p_{*}q^{*}(H^{n}_{prim} (X))\\ & =p_{*}(p^{*}{\iota}^{*}H^{3(n-2)} (G_{1,n+1})\cdot q^{*}H^{n}_{prim} (X))\\ & =p_{*}(\overline{\iota}^{*}\overline{p}^{*}H^{3(n-2)} (G_{1,n+1})\cdot q^{*}H^{n}_{prim} (X)) \end{split}\]
where the second equality uses the projection formula and where the third equality uses the commutativity of the left square in (\ref{comm}). Thus, (\ref{goal0}) reduces to showing that 
\begin{equation}\label{dot0}
p_{*}(\overline{\iota}^{*}H^{3(n-2)} (\pit_{G_{1,n+1}}(\E))\cdot q^{*}H^{n}_{prim} (X)) = 0
\end{equation}
For this, we observe that $\pit_{G_{1,n+1}}(\E) \xrightarrow{q} \pit^{n+1}$ is the projective bundle $\pit(T\pit^{n+1}) \to \pit^{n+1}$, where $T\pit^{n+1}$ is the tangent bundle of $\pit^{n+1}$. To see this, note that for $x \in X(\cit)$ 
\[ q^{-1} (x) = \{ (y, \ell) \in \pit^{n+1} \times G_{1,n+1} \ | \ x, y \in \ell \} \]
which one realizes as the fiber over $x$ of the projective bundle $\pit(T\pit^{n+1}) \to \pit^{n+1}$ using the exact sequence:
\[ 0 \to \mathcal{O}_{\pit^{n+1}} \to \mathcal{O}_{\pit^{n+1}}(1)^{\oplus n+2} \to T\pit^{n+1} \to 0 \]
Then, we have the diagram:
\[
\begin{tikzcd}
\pit(\E) \arrow{r}{q} \arrow{d}{\iota} & X \arrow{d}{j} \\
\pit_{G_{1,n+1}}(\E) \cong \pit(T\pit^{n+1}) \arrow{r}{\overline{q}} & \pit^{n+1}
\end{tikzcd}\]
To show that 
\begin{equation}\label{cdot0} \overline{\iota}^{*}H^{3(n-2)} (\pit_{G_{1,n+1}}(\E))\cdot q^{*}H^{n}_{prim} (X) = 0
\end{equation}
we observe that 
\[ H^{3(n-2)} (\pit_{G_{1,n+1}}(\E)) = H^{3(n-2)} (\pit(T\pit^{n+1})) = \bigoplus_{k=0}^{n+1} {\overline{q}}^{*}H^{3(n-2) - 2k} (\pit^{n+1})\cdot c_{1}(\mathcal{O}_{T\pit^{n+1}}(1)))^{k} \]
Then, for each $k$, we have 
\[\begin{split} {\overline{\iota}}^{*}({\overline{q}}^{*}H^{3(n-2) - 2j} (\pit^{n+1})\cdot c_{1}(\mathcal{O}_{T\pit^{n+1}}(1))^{k})\cdot{q}^{*}(H^{n}_{prim} (X)) \\ = {q}^{*}(j^{*}H^{3(n-2) - 2j} (\pit^{n+1})\cdot H^{n}_{prim} (X))\cdot {\overline{\iota}}^{*}c_{1}(\mathcal{O}_{{T\pit^{n+1}}}(1))^{k} = 0 \end{split} \]
since $H^{n}_{prim} (X)\cdot H = 0$. This gives (\ref{cdot0}) and, hence, (\ref{goal0}).
\end{proof}
\end{Lem}

\noindent The next lemma gives a characterization of the cohomology of $F$ except in the middle degree. We first introduce the following notation. 
\begin{Not} Denote by $(-,-)_{k}: H^{k} (F) \otimes H^{k} (F) \to \qit$ the pairing:
\[(\alpha, \alpha')_{k} := \alpha\cdot \alpha'\cdot c_{1}(\E)^{2(n-2) - k}\]
for $\alpha, \alpha' \in H^{k} (F)$. By Poincar\'e duality and the Hard Lefschetz decomposition, $(-,-)_{k}$ is non-degenerate.
\end{Not} 
\begin{Lem}\label{cor-lem} For $0 \leq k \neq 4(n-2)$ with $k \neq 2(n-2)$,
\[H^{k}(F) :=   \left\{\begin{array}{cl} 
 L^{s}\Gamma_{*}H^{n}_{prim} (X) & \text{ for $n$ odd and $k = n-2+2s$, $0\leq s \leq n-2$}\\
L^{s}\Gamma_{*}H^{n}_{prim} (X) \oplus \iota^{*}H^{k} (G_{1,n+1}) & \text{ for $n$ even and $k = n-2+2s$, $0\leq s \leq n-2$}\\
\iota^{*}H^{k} (G_{1,n+1}) & \text{ otherwise } \end{array} \right. \]
where $L$ denotes the Lefschetz operator for $c_{1}(\E)$ and where $L^{s}\Gamma_{*}H^{n}_{prim} (X) \oplus \iota^{*}H^{k} (G_{1,n+1})$ is an orthogonal decomposition with respect to $(-,-)_{k}$.
\begin{proof} For the first, note that, by the Hard Lefschetz theorem,
\begin{equation} L^{s}: H^{n-2} (F) \to H^{k} (F)\label{lef} \end{equation}
is injective. Then, using (\ref{cohom-1}) it follows that $H^{n-2} (F)$ and $H^{k} (F)$ have the same rank. The third follows from Corollary \ref{cor-cor}. For the second, first note that $\iota^{*}H^{k} (G_{1,n+1})$ and $L^{s}\Gamma_{*}H^{n}_{prim}(X)$ are orthogonal for $(-,-)_{k}$. Indeed, let $\alpha \in L^{s}\Gamma_{*}H^{n}_{prim}(X)$ and write $\alpha = \alpha'\cdot c_{1}(\E)^{s}$ for $\alpha' \in \Gamma_{*}H^{n}_{prim}(X)$ and let $\beta \in \iota^{*}H^{k} (G_{1,n+1})$. Then, $\beta\cdot c_{1}(\E)^{s} \in \iota^{*}H^{3(n-2)} (G_{1,n+1})$ and by (\ref{goal}) the subspaces $\iota^{*}H^{3(n-2)} (G_{1,n+1})$ and $\Gamma_{*}H^{n}_{prim}(X)$ are orthogonal. It follows that 
\[ (\alpha, \beta)_{k} = \alpha'\cdot(\beta\cdot c_{1}(\E)^{s}) = 0 \]
Now, from the injectivity of (\ref{lef}), it follows that  
\[\text{rank}(L^{s}\Gamma_{*}H^{n}_{prim}(X)) = \text{rank}(\Gamma_{*}H^{n}_{prim}(X)) = \text{rank}(H^{n}_{prim}(X)) \]
By Corollary \ref{first}, we have $\text{rank}(\iota^{*}H^{k} (G_{1,n+1})) \geq a_{\frac{k}{2}}$. Moreover, using (\ref{cohom-1}) we obtain
\[ \text{rank}(H^{k} (F)) =  \text{rank}(H^{n}_{prim}(X))+ a_{\frac{k}{2}}\]
Thus, it follows that
\begin{equation} \text{rank}(\iota^{*}H^{k} (G_{1,n+1})) + \text{rank}(L^{s}\Gamma_{*}H^{n}_{prim}(X)) \geq \text{rank}(H^{k} (F)) \label{ineq}\end{equation}
Since $\iota^{*}H^{k} (G_{1,n+1})$ and $L^{s}\Gamma_{*}H^{n}_{prim}(X)$ are orthogonal with respect to $(-,-)_{k}$, we deduce that (\ref{ineq}) is an equality and, hence, that 
\[H^{k} (F) = L^{s}\Gamma_{*}H^{n}_{prim} (X) \oplus \iota^{*}H^{k} (G_{1,n+1})\]
\end{proof}
\end{Lem}

\begin{Cor}\label{last-cor} For all $0 \leq k \leq 4(n-2)$ with $k \neq 2(n-2)$, \[\text{rank}(\iota^{*}H^{k} (G_{1,n+1})) = a_{\frac{k}{2}} = a_{\frac{4(n-2)-k}{2}} =\text{rank}(\iota^{*}H^{4(n-2)-k} (G_{1,n+1}))\]  Moreover, the pairing:
\[ \iota^{*}H^{k} (G_{1,n+1}) \otimes \iota^{*}H^{4(n-2)-k} (G_{1,n+1}) \to \qit \]
induced by the cup product is non-degenerate.
\begin{proof} For the first statement, note by Corollary \ref{first} we have 
\[\text{rank}(\iota^{*}H^{k} (G_{1,n+1})) \geq a_{\frac{k}{2}} \] However, the last few sentences of the proof of Lemma \ref{cor-lem} show that this must be an equality. It suffices to prove the second statement in the case that $k \leq 2(n-2)$. To this end, suppose that $\alpha \neq 0 \in H^{k} (G_{1,n+1})$. The proof of Corollary \ref{first} then shows that
\[\alpha \cdot F \neq 0 \in H^{k+8} (Gr(1, n+1))\] 
So, there exists $\beta \neq 0 \in H^{4n-(k+8)} (Gr(1, n+1))$ such that $\alpha\cdot \beta \cdot F \neq 0$. Thus, by the projection formula,
\[ \iota_{*}(\iota^{*}\alpha\cdot\iota^{*}\beta) = \alpha\cdot \beta \cdot F \neq 0 \in H^{4n}(Gr(1, n+1))  \]
from which it follows that $\iota^{*}\alpha\cdot\iota^{*}\beta \neq 0$.
\end{proof}
\end{Cor}

\section{Proof of Theorem \ref{main}}
\noindent Since we are working over $\cit$, it will suffice to prove that for $0 \leq k \leq 2(n-2)$ there exist correspondences $\Gamma_{k} \in CH^{k}(F \times F) = Cor^{k-2(n-2)} (F, F)$ for which
\[\Gamma_{k*}: H^{4(n-2)-k} (F) \to H^{k} (F) \]
is an isomorphism. Indeed, according to Theorem 2.9 of \cite{Kl}, this implies both the Lefschetz standard and the K\"unneth conjectures. Then, \cite{Gr} shows that the Lefschetz standard conjecture implies that $\sim_{num} = \sim_{hom}$. What remains then is to construct the required correspondences $\Gamma_{k}$. To this end, let $L \in CH^{2(n-2)+1} (F \times F)=Cor^{1} (F, F)$ denote the Lefschetz correspondence for $c_{1} (\E)$. Also, let 
\[ \delta_{prim} := \Delta_{X} - \sum_{0 \leq r \leq n} \frac{1}{3} H^{r} \times H^{n-r} \in CH^{n} (X \times X) =Cor^{0} (X, X) \]
be the correspondence for which $\delta_{prim*}H^{*}(X) = H^{n}_{prim}(X)$. Then, we can define the following correspondence in $CH^{k} (F \times F)$:
\begin{equation} \Gamma_{k}' :=   \left\{\begin{array}{cl} L^{s} \circ\Gamma\circ\delta_{prim} \circ \prescript{t}{}{\Gamma}\circ L^{s} & \text{ for $k = n-2+2s, 0 \leq s \leq [\frac{n-2}{2}]$ }\\
0 & \text{ for all other $k$ } \end{array} \right. \label{gamma'}\end{equation}
\begin{Lem}\label{gamma'lem} When $\Gamma_{k}' \neq 0$, we have 
\begin{enumerate}
\item\label{gamma'1} $\Gamma_{k*}'(H^{k} (F)) = L^{s}\Gamma_{*}H^{n}_{prim} (X)$.  
\item\label{gamma'2} $\Gamma_{k*}'(\iota^{*}H^{4(n-2)-k} (G_{1,n+1})) =0 $
\end{enumerate}
\begin{proof} The statement of (\ref{gamma'1}) will follow from the fact that \begin{equation}\delta_{prim} \circ \prescript{t}{}{\Gamma}\circ L^{s}: H^{4(n-2)-k} (F) \to H^{n}_{prim}(X)\label{red'} \end{equation} is surjective for $k =n-2+2s$ (with $0 \leq s \leq [\frac{n-2}{2}]$). To this last end, note that the hard Lefschetz theorem implies that $L^{s}H^{4(n-2)-k}(F) = H^{3(n-2)}(F)$; i.e., 
\[ L^{s}: H^{4(n-2)-k} (F) \to H^{3(n-2)}(F)\] 
is surjective. Moreover, $\delta_{prim} \circ \prescript{t}{}{\Gamma} : H^{3(n-2)} (F) \to H^{n}_{prim}(X)$ is surjectve since its Poincar\'e dual \[\prescript{t}{}{(\delta_{prim} \circ \prescript{t}{}{\Gamma})} = \Gamma\circ\delta_{prim}: H^{n}_{prim} (X) \to H^{n-2}(F)\] is injective by Lemma \ref{second}. For the statement of (\ref{gamma'2}), it will suffice to show that 
\begin{equation} \iota^{*}H^{4(n-2)-k} (G_{1,n+1}) \subset \ker\{\delta_{prim} \circ \prescript{t}{}{\Gamma}\circ L^{s}: H^{4(n-2)-k} (F) \to H^{n}_{prim}(X)\}\label{goal}\end{equation}
or equivalently that 
\[ \iota^{*}H^{3(n-2)} (G_{1,n+1}) \subset \ker\{\delta_{prim} \circ \prescript{t}{}{\Gamma}: H^{3(n-2)} (F) \to H^{n}_{prim}(X)\}\]
For this, note that
\[ \ker\{\delta_{prim} \circ \prescript{t}{}{\Gamma}: H^{3(n-2)} (F) \to H^{n}_{prim}(X)\}) = (\Gamma_{*}H^{n}_{prim}(X))^{\perp} \]
where $(\phantom{\text{bl}})^{\perp}$ is with respect to the cup product. The statement of (\ref{goal}) then follows from the last statement of Lemma \ref{second}. 
\end{proof}
\end{Lem} 
\begin{Lem}\label{gamma''lem} For $0 \leq k \leq 2(n-2)$, there exists a correspondence $\Gamma_{k}'' \in CH^{k} (F \times F)$ such that:
\begin{enumerate} 
\item\label{gamma''1} $\Gamma_{k*}''(\iota^{*}H^{4(n-2)-k} (G_{1,n+1})) = \iota^{*}H^{k} (G_{1,n+1})$
\item\label{gamma''2} If $n$ is even and $k = n-2 +2s$ for some $s$, $\Gamma_{k*}''(L^{n-2-s}\Gamma_{*}H^{n}_{prim} (X)) = 0$
\end{enumerate}

\begin{proof}For $k$ odd, set $\Gamma_{k}'' =0$. For $k$ even, we let \[\{ \alpha_{k, m} \} \subset \iota^{*}H^{k} (G_{1,n+1}) \] be a basis, for which there are lifts $\tilde{\alpha}_{k,m} \in \iota^{*}CH^{\frac{k}{2}} (G_{1,n+1})$. Then, we can set 
\[ \Gamma_{k}'' = \sum_{m=1}^{[\frac{k+4}{4}]} \tilde{\alpha}_{k, m} \times \tilde{\alpha}_{k, m} \in CH^{k} (F \times F)= Cor^{k-2(n-2)} (F, F) \]
By Corollary \ref{last-cor}, there exists a dual basis 
\[\{ \beta_{k, m} \} \subset \iota^{*}H^{4(n-2)-k} (G_{1,n+1}) \]
with respect to the cup product. A standard computation then shows that
\[ \Gamma_{k*}''(\beta_{k, n}) = \delta_{m,n} \cdot \alpha_{k,m} \in \iota^{*}H^{k} (G_{1,n+1}) \]
where $\delta_{m,n} = 1$ if $m=n$ and $0$ otherwise. This proves (\ref{gamma''1}). For (\ref{gamma''2}), we note by Lemma \ref{cor-lem} that $\iota^{*}H^{k} (G_{1,n+1})$ and $L^{n-2-s}\Gamma_{*}H^{n}_{prim} (X)$ are orthogonal with respect to the cup product (since $\iota^{*}H^{k} (G_{1,n+1})$ and $L^{s}\Gamma_{*}H^{n}_{prim} (X)$ are orthogonal with respect to $(-,-)_{k}$). Thus, for all $m$ and all $\gamma \in L^{s}\Gamma_{*}H^{n}_{prim} (X)$, we have $\alpha_{k,m}\cdot \gamma =0$, from which it follows that $\Gamma_{k*}''(\gamma) =0$.
\end{proof}
\end{Lem}
\noindent To complete the proof, we set
\begin{equation} \Gamma_{k} = \Gamma_{k}' + \Gamma_{k}'' \in CH^{k} (F \times F)= Cor^{k-2(n-2)} (F, F)
\end{equation}
Then, we have the following lemma:
\begin{Lem} $\Gamma_{k}: H^{4(n-2)-k} (F) \to H^{k} (F)$ is an isomorphism for $0 \leq k \leq 2(n-2)$.
\begin{proof} It suffices to show that $\Gamma_{k}$ is surjective (since the ranks of $H^{4(n-2)-k} (F)$ and $H^{k} (F)$ are the same). To this end, we observe from Lemma \ref{cor-lem}:
\begin{equation} H^{k}(F) :=   \left\{\begin{array}{cl} L^{s}\Gamma_{*}H^{n}_{prim} (X) \oplus \iota^{*}H^{k} (G_{1,n+1}) & \text{ for $n$ even and $k = n-2+2s$}\\
 L^{s}\Gamma_{*}H^{n}_{prim} (X) & \text{ for $n$ odd and $k = n-2+2s$}\\
\iota^{*}H^{k} (G_{1,n+1}) & \text{ otherwise } \end{array} \right. \end{equation}
When $H^{4(n-2)-k}(F) = \iota^{*}H^{4(n-2)-k} (G_{1,n+1})$, surjectivity follows from the fact that $\Gamma_{k}' = 0$ in this case and Lemma \ref{gamma''lem}. When $n$ is odd and $k = n-2+2s$, surjectivity follows from the fact that $\Gamma_{k}'' = 0$ and Lemma \ref{gamma'lem}. What remains then is to prove surjectivity in the case that $n$ is even and $k = n-2+2s$. For this, we define
\[ \begin{split}  V' = L^{n-2-s}\Gamma_{*}H^{n}_{prim} (X), & \ V'' = \iota^{*}H^{4(n-2)-k} (G_{1,n+1}),\\ W' = L^{s}\Gamma_{*}H^{n}_{prim} (X),  & \   W'' = \iota^{*}H^{k} (G_{1,n+1}) \end{split}\]
Then, we observe that \[ \Gamma_{k*}'(V')= W', \ \Gamma_{k*}'(V'') = 0, \ \Gamma_{k*}''(V') = 0, \ \Gamma_{k*}''(V'') = W''\] Since $H^{4(n-2)-k}(F) = V' \oplus V''$ and $H^{k}(F) = W' \oplus W''$, the lemma then follows from the following (essentially trivial) fact from linear algebra:
\begin{claim} Let $V= V' \oplus V''$ and $W= W' \oplus W''$ are vector spaces. Suppose that $T', T'': V \to W$ are linear maps for which $T' (V') = W'$ , $T' (V'') = 0$, $T''(V') = 0$ and $T''(V'') = W''$.
 Then, $T = T' + T'': V \to W$ is surjective.
\end{claim}
\end{proof}
\end{Lem}

\end{document}